\documentclass{amsart}
\newtheorem{theorem}{Theorem}[section]
\newtheorem{lemma}[theorem]{Lemma}
\newtheorem{prop}[theorem]{Proposition}
\newtheorem{rem}[theorem]{Remark}
\newtheorem{coro}[theorem]{Corrolary}
\newtheorem{defi}[theorem]{Definition}
\newcommand{\bproof}{{\bf Proof:~~}}
\newcommand{\eproof}{{\vrule height8pt width5pt depth0pt}\vspace{3mm}}

\newcommand{\mref}[1]{{(\ref{#1})}}
\newcommand{\reflemm}[1]{Lemma~\ref{#1}}
\newcommand{\refrem}[1]{Remark~\ref{#1}}
\newcommand{\reftheo}[1]{Theorem~\ref{#1}}

\newcommand{\refcoro}[1]{Corollary~\ref{#1}}
\newcommand{\refprop}[1]{Proposition~\ref{#1}}
\newcommand{\refsec}[1]{Section~\ref{#1}}

\newcommand{\mI}{\mathcal{I}}
\newcommand{\mT}{\mathcal{T}}
\newcommand{\mR}{\mathcal{R}}

\newcommand{\mM}{\mathcal{M}}
\newcommand{\mP}{\mathcal{P}}

\newcommand{\mW}{\mathcal{W}}
\newcommand{\mB}{\mathcal{B}}

\newcommand{\llg}{\lambda}
\newcommand{\ag}{\alpha}
\newcommand{\bg}{\beta}
\newcommand{\sg}{\sigma}

\newcommand{\og}{\omega}
\newcommand{\Og}{\Omega}
\newcommand{\pdh}{\partial}
\newcommand{\RR}{{\rm I\kern -1.6pt{\rm R}}}

\title[A family of inversion formulas in Thermoacoustic Tomography]
      {A family of inversion formulas in thermoacoustic tomography}

\author[Linh V. Nguyen]{}

\subjclass{Primary: 92C55, 35L05, 35R30; Secondary: 44A99, 45Q05}
 \keywords{Tomography, Inversion formulas, Spherical means, Wave equation}

\begin{document}

\maketitle

%% Enter the first author's name and address:
\centerline{\scshape Linh V. Nguyen }
\medskip
{\footnotesize
 %% please put the address of the first author
 \centerline{Department of Mathematics, Texas A \& M University }
   \centerline{Mailstop 3368, College Station, TX 77843-3368, USA}
} %% Do not forget to end the {\footnotesize by the sign }

\medskip

\bigskip

\begin{abstract} We present a family of closed form inversion formulas in thermoacoustic tomography in the case of a constant sound speed. The formulas are presented in both time-domain and frequency-domain versions. As special cases, they imply most of the previously known filtered backprojection type formulas.
\end{abstract}

%\tableofcontents

\section{Introduction and statement of main results}\label{intro}
Thermoacoustic tomography (TAT) has recently attracted considerable attention as a promising method of biomedical imaging \cite{XW06,WWu,Ora02,Ora03,KLFA}. We give here a brief introduction to the mathematical model of TAT.

An electromagnetic (EM) pulse in visible light or radiofrequency range is sent through a biological object of interest. A fraction of the EM energy is absorbed by the tissues throughout the object. It is known that cancerous tissues absorb much more EM energy than healthy ones (e.g., \cite{XW06,KLFA}). Thus, knowing the absorption distribution function $a(x)$ would be useful for detecting cancerous tissues. The EM energy absorbed by the tissue causes its thermoelastic expansion, which in turn releases a pressure (ultrasound) wave $u(x,t)$ propagating through the body. This pressure is then measured by transducers located on an {\bf observation surface} $S$ surrounding the object. The initial pressure $f(x) = u(x,0)$ is roughly proportional to $a(x)$. One now concentrates on the recovery of $f(x)$ from the measured data $g:=u|_{S \times [0,\infty)}$. There has been a number of papers devoted to this reconstruction (see reviews in \cite{FR07,FR09,KKun,AKKun}). All known closed form inversion formulas assume that the object is acoustically homogeneous, i.e. the sound speed is constant. We will assume that, after choosing appropriate units, the sound speed is equal to one.

The standard mathematical model of TAT is \cite{XW02,DSK,Tam}
\begin{eqnarray}\label{wave}
\left \{\begin{array}{l} u_{tt}(x,t) - \bigtriangleup  u(x,t) =0,~ x \in \RR^n,~ t \geq 0, \\ u(x,0) =f(x), \quad u_t(x,0) = 0,\\u(y,t) =g(y,t), \mbox{ for } y \in S,~ t \geq 0. \end{array} \right.
\end{eqnarray}
Here $g(x,t)$ is the measured data and $f$ is the function to be reconstructed. In other words, one needs to invert the operator $\mT: f \mapsto g$. Since $\mT$ is known to be invertible from the left only, different inversion formulas exist. In this text, we develop a family of explicit closed form inversion formulas in the case when $S$ is the unit sphere centered at the origin and $f \in C_0^\infty(\overline{B})$, where $B$ is the open unit ball enclosed by $S$.

We will also need to deal with some other operators closely related to $\mT$. Consider the wave equation \begin{eqnarray}\label{wave1} \left \{\begin{array}{l} u_{tt}(x,t) - \bigtriangleup  u(x,t) =0,\quad  x \in \RR^n,\quad t \geq 0, \\ u(x,0) =0, \quad u_t(x,0) = f(x),\\u(y,t) =g(y,t), \mbox{ for } y \in S, t \geq 0. \end{array} \right.\end{eqnarray} One can now define the the operator $\mP$ that maps function $f$ into $g$: $\mP(f)=g$. Then $\mP$ is related to $\mT$ as follows: $$\mP(f)(y,t) = \int\limits_0\limits^t \mT(f)(y,\tau) d\tau \mbox{ and }\mT(f)(t)=\pdh_t \mP(f)(t).$$
We also introduce the spherical Radon transform $\mR_S$ with centers on $S$ by $$\mR_S(f)(y,t) = \int\limits_{S^{n-1}}f(y+r \og) r^{n-1} d\sg(\og),\mbox{ for all } (y,t) \in S \times \RR_+.$$ Finally, $\mM_S$ is the spherical means operator with centers on $S$:
$$\mM_S(f)(y,t) =\frac{1}{\og_n} \int\limits_{S^{n-1}}f(y+r \og)d\sg(\og), \mbox{ for all } (y,t) \in S \times \RR_+.
$$
In these formulas, $d\sg(\og)$ is the standard surface measure on the unit sphere $S^{n-1}\subset \RR^n$ and $\og_n$ is the total measure of the sphere.

In various works cited below, different operators from this list were considered. However, due to known explicit connections between $\mT$ and $\mP$, $\mR_S$, and $\mM_S$ (e.g., \cite{CH2,Jo}), inversion formulas for these operators are closely related.

The first such formulas, in odd dimensions, were obtained by Finch, Patch and Rakesh \cite{FPR} using some trace identities for the wave operator. Xu and Wang \cite{XW05} derived a different formula for $n=3$ by working in the frequency domain. A formula for all dimensions, which coincides with that of \cite{XW05} when $n=3$, was presented by Kunyansky \cite{Kun07}. Its derivation is based upon some symmetry relation for special functions. By applying the same method as in \cite{FPR}, Finch, Haltmeier and Rakesh \cite{FHR} obtained inversion formulas for even dimensions, which involves the data measured for an infinite time period. The authors of \cite{FHR} also derived another type of inversion formulas for even $n$, which uses only the data measured for a finite period of time (which we will refer as "finite-time formulas for even dimensions").

In spite of availability of several types of closed form inversion formulas, some questions have remained unanswered. For instance, it was not clear, what is the relation, if any, between formulas of \cite{FHR,FPR} and \cite{Kun07,XW05}, (which are known to be not equivalent outside the range of the operator $\mT$). The same applies to the two types of inversion formulas derived in \cite{FHR} for even dimensions.

The goal of this work is to obtain a unified family of closed form inversion formulas, which in particular would produce formulas of \cite{FPR,FHR,Kun07,XW05} as particular cases.

We now formulate the main results of the paper. Consider the function \begin{equation} \label{green} G(s,\llg) = \frac{i}{4} \left( \frac{\llg}{2 \pi s}\right)^{\frac{n-2}{2}} H^{(1)}_{\frac{n-2}{2}}(\llg s),\end{equation} where $H^{(1)}_{\frac{n-2}{2}}$ is the Hankel function of the first kind. Then (e.g., \cite{AGM}) $\Phi(x,y,\llg) =G(|x-y|,\llg)$ is the solution for the Helmholtz equation \begin{equation}\label{he}\Delta U(y) + \llg^2 U(y) =-\delta(y-x),\end{equation} obtained by limiting absorption. For odd $n$ and any $s>0$, $G(s,.) \in C^\infty(\RR)$. For even $n$, $G(s,.) \in C^\infty(\RR \setminus \{0\})$ with a logarithmic singularity at zero \cite{AGM}.

Let $g(y,t)$ be the function in (\ref{wave}) that represents the thermoacoustic data and $g_0(y,t)$ be its even extension with respect to $t$. Let also $\hat{g}_0$ be the Fourier transform of $g_0$ with respect to time:
$$
\hat{g}_0(y,\llg) = \frac{1}{2 \pi}\int\limits_{\RR}g_0(y,t)e^{i \llg t} dt.
$$
Then, our family of inversion formula reads as follows:
\begin{theorem} \label{geninf}  Suppose that $f \in C^\infty_0(\overline{B})$ and $g=\mT(f)$. Then for any $x \in B$ and $\xi \in \RR^n$, the following equality holds \begin{eqnarray} \label{genform} f(x)&=& - 2 \int\limits_{S} \left.\left(\frac{d}{ds} \int\limits_{\RR} \overline{G}(s,\llg)\hat{g}_0(y,\llg) d\llg \right)\right|_{s=|x-y|} \frac{\left<y-x,y-\xi \right>}{|x-y|} d\sg(y),\end{eqnarray} where $\overline{G}$ is the complex conjugate of $G$.\end{theorem}

One can rewrite this inversion formula without going to the frequency domain. Namely, we define a transform $\mW$ of a function $v \in C^\infty[0,\infty)$ as follows (as long as the expressions involved make sense): \begin{eqnarray}\label{mw} \mW(v)(s) :=\left\{\begin{array}{l} c_n \left(\frac{1}{s} \frac{d}{ds} \right)^{\frac{n-2}{2}} \int\limits_{s}^\infty \frac{v(t)}{\sqrt{t^2-s^2}}  dt, \mbox{ if $n$ is even}, \\ c_n \left(\frac{1}{s} \frac{d}{ds} \right)^{\frac{n-3}{2}} \left(\frac{v(s)}{s}\right), \mbox{ if $n$ is odd}.
\end{array} \right.\end{eqnarray} Here \begin{eqnarray} \label{cn} c_n := \left\{\begin{array}{l} \frac{(-1)^{\frac{n-2}{2}}}{(2 \pi)^{\frac{n}{2}}},\mbox{ if $n$ is even}, \\ \frac{(-1)^{\frac{n-3}{2}}}{2(2\pi) ^{\frac{n-1}{2}}},\mbox{ if $n$ is odd}.\end{array} \right.\end{eqnarray}
 (As noted in \refsec{remm}, $\mW$ is a known transform that intertwines the second derivative and the Bessel operator.) One can now obtain another representation of the inner integral in \mref{genform}: $$\int\limits_{\RR} \overline{G}(s,\llg)\hat{g}_0(y,\llg) d\llg= \mW(g)(y,s),$$ where $\mW(g)(y,s):=\mW(g_y)(s)$ with $g_y(t) :=g(y,t)$. Then \reftheo{geninf} is equivalent to
\begin{theorem} \label{genin} Suppose that $f \in C^\infty_0(\overline{B})$ and $g=\mT(f)$. Then for any $x \in B$ and $\xi \in \RR^n$, the following equality holds \begin{eqnarray} \label{invf} f(x) &=&- 2 \int\limits_{S} \left.\left(\frac{d}{ds} \mW(g)(y,s) \right)\right|_{s=|x-y|} \frac{\left<y-x,y-\xi \right>}{|x-y|} d\sg(y).\end{eqnarray}
\end{theorem}

One can look at the inner integral in \mref{genform} in a different manner yet. Namely, due to the relation between $\mT$ and $\mR_S$, it can be shown that  $$\int\limits_{\RR} \overline{G}(s,\llg)\hat{g}_0(y,\llg) d\llg =\frac{2}{\pi} \int\limits_0\limits^\infty \llg \mR(s,\llg) \int\limits_0\limits^\infty \mR_S(f)(y,r) \mI(r,\llg) dr,$$ where $\mR$ and $\mI$ are the real and imaginary parts of $G$.
Then \reftheo{geninf} implies
\begin{theorem} \label{Kunt} Suppose that $f \in C^\infty_0(\overline{B})$ and $g=\mR_S(f)$. Then for any $x \in B$ and $\xi \in \RR^n$, the following equality holds \begin{eqnarray}\label{knintro} f(x)= - \frac{4}{\pi} \int\limits_{S} \left(\frac{d}{ds} K_n(y,s) \right)_{s=|x-y|} \frac{\left<y-x,y-\xi\right>}{|x-y|} d\sg(y),\end{eqnarray} where $$K_n(y,s) =\int\limits_0\limits^\infty \llg \mR(s,\llg) \int\limits_0\limits^\infty g(y,r) \mI(r,\llg) dr d\llg.$$
\end{theorem}

The key ingredient in the proof of Theorem \ref{geninf} is the following identity, which is equivalent to a range description of operator $\mT$ (see \refrem{rangerem}): \begin{theorem} \label{main} Suppose that $\Og$ is a ball, $f \in C^\infty_0(\overline{\Og})$, and $u$ solves (\ref{wave}). Let $u_0$ be the even extension of $u$ with respect to $t$. Then for all $x \in \Og$ we have \begin{eqnarray} \label{mainr} \int\limits_{\pdh \Og} \int\limits_{\RR} \overline{G}(|x-y|,\llg)\hat{u}_0(y,\llg) d\llg d\sg(y)=0.\end{eqnarray} \end{theorem}

This statement implies that adding to any inversion formula an expression
\begin{eqnarray*}
A\int\limits_{S} \int\limits_{\RR} \overline{G}(|x-y|,\llg)\hat{g}_0(y,\llg) d\llg d\sg(y)
\end{eqnarray*}
with an arbitrary operator $A$ also produces an inversion formula (since the added term vanishes on the range of the operator to be inverted). In particular, if $A$ is the operator of multiplication by an arbitrary function $\varphi(x)$, one can add the expression
\begin{eqnarray*} \varphi(x)\int\limits_{S} \int\limits_{\RR} \overline{G}(|x-y|,\llg)\hat{g}_0(y,\llg) d\llg d\sg(y).\end{eqnarray*} to \mref{invf} and \mref{knintro} to obtain the following inversion formulas:

\begin{coro} \label{fam}  Suppose that $f \in C^\infty_0(\overline{B})$ and $g=\mT(f)$. Let $\varphi$ be an arbitrary function defined on $B$. Then for any $x \in B$ and $\xi \in \RR^n$, the following equalities hold \begin{eqnarray}\label{fam1} \nonumber f(x) &=& - \frac{4}{\pi} \int\limits_{S} \left(\frac{d}{ds} K_n(y,s) \right)_{s=|x-y|} \frac{\left<y-x,y-\xi\right>}{|x-y|} d\sg(y),\\ &+& \varphi(x) \int\limits_{S} K_n(y,|x-y|) d\sg(y),\\ \label{fam2} \nonumber f(x) &=&- 2 \int\limits_{S} \left.\left(\frac{d}{ds} \mW(g)(y,s) \right)\right|_{s=|x-y|} \frac{\left<y-x,y-\xi \right>}{|x-y|} d\sg(y) \\ &+& \varphi(x) \int\limits_{S} \mW(g)(y,|x-y|) d\sg(y). \end{eqnarray}\end{coro}

Now, by suitable choices of $\xi$ and $\varphi$ in \mref{fam1} and \mref{fam2}, one can recover the inversion formulas known in the literature.
If we let $\xi=x$, and $\varphi=-2(n-2), -2(n-1)$, or $-2n$ in \mref{fam2}, we obtain correspondingly the formulas derived in \cite[Theorem 4]{FHR} and \cite[Theorem 3]{FPR}:
\begin{prop}\label{fprop} Let $f \in C_0^\infty(\overline{B})$. Then for any $x \in B$:
\begin{eqnarray} \label{f1} f(x) &=&-2\left(\mP^*t\pdh_t^2 \mP f\right)(x)\\ \label{f2} f(x) &=&  -2\left(\mP^* \pdh_t t \pdh_t \mP f \right)(x)\\ \label{f3} f(x) &=&  -2\left(\mP^* \pdh_t^2 t \mP f \right)(x)\end{eqnarray}
Here $\pdh_t= \frac{d}{dt}$ is the derivative with respect to $t$, and $\mP^*$ is the $L^2-$adjoint of $\mP$.
\end{prop}

Choosing $\xi=x$ and $\varphi=0$ in \mref{fam1}, one gets a finite-time inversion formula for even dimensions similar to the second one in \cite[Theorem 2]{FHR}:
\begin{prop}\label{ffinite} Assume that $n$ is even, $f \in C_0^\infty(\overline{B})$, and $g=\mM_S (f)$. Then for all $x\in B$, $f(x)$ is equal to \begin{eqnarray*} \frac{(-1)^{\frac{n-2}{2}}\og_n}{(2\pi)^{n}} \int\limits_{S} \int\limits_0\limits^2 \left[\pdh_r r\left(\pdh_r \frac{1}{r}\right)^{n-1}r^{n-1}g\right](y,r) \ln \left|r^2-|x-y|^2\right| dr d\sg(y),\end{eqnarray*} where $\og_n$ is the surface area of the unit sphere $S^{n-1}$.\end{prop}

In particular, when $n=2$, we obtain the second formula in \cite[Theorem 1]{FHR}:
\begin{coro}Assume that $n=2$, $f \in C^\infty_0(\overline{B})$, and $g=\mM_S(f)$. Then for any $x \in B$, \begin{eqnarray}\label{n2} f(x) =  \frac{1}{2 \pi} \int\limits_{S} \int\limits_0\limits^2 \left[\pdh_r r \pdh_r g\right](y,r) \ln \left|r^2 -|x-y|^2\right| dr d\sg(y). \end{eqnarray} \end{coro}

Finally, let $$J(s)=\frac{J_{\frac{n-2}{2}}(s)}{s^{\frac{n-2}{2}}},~~N(s)=\frac{N_{\frac{n-2}{2}}(s)}{s^{\frac{n-2}{2}}},$$ where $J_{\frac{n-2}{2}}$ and $N_{\frac{n-2}{2}}$ are the Bessel and Neumann functions of order $\frac{n-2}{2}$. Choosing $\xi=0$ and $\varphi=0$ in \mref{fam1}, one arrives at:
\begin{prop}\label{kuneq} Assume that $f\in C_0^\infty(\overline{B})$ and $g=\mR_S(f)$. Then for any $x \in B$, \begin{eqnarray}\label{mine} f(x) = \frac{- 1}{2 (2\pi)^{n-1}} \nabla_x \int\limits_{S} n(y) k_n(y,|x-y|) d\sg(y).\end{eqnarray} where
\begin{eqnarray*}k_n(y,s) &=&  \int\limits_0\limits^\infty \llg^{2n-3} N(s\llg) \int\limits_0\limits^2 g(y,r) J(r \llg) dr d\llg.\end{eqnarray*}  \end{prop}
This inversion formula is equivalent to the one obtained in \cite{Kun07}. Indeed, let  \begin{eqnarray*} h(y,s) &=& \int\limits_0\limits^\infty \llg^{2n-3} N(s\llg) \int\limits_0\limits^2 g(y,r) J(r \llg) dr d\llg \\ &-& \int\limits_0\limits^\infty \llg^{2n-3} J(s\llg) \int\limits_0\limits^2 g(y,r) N(r \llg) dr d\llg. \end{eqnarray*}
It can be shown that $h(y,s) = 2 k_n(y,s)$, and thus \refprop{kuneq} implies the following result of \cite{Kun07}:
\begin{prop}\label{kunprop} \cite{Kun07} Assume that $f\in C_0^\infty(\overline{B})$ and $g=\mR_S(f)$. Then for any $x \in B$, \begin{eqnarray}\label{kunintro} f(x) = \frac{- 1}{4 (2\pi)^{n-1}} \nabla_x \int\limits_{S} n(y) h(y,|x-y|) d\sg(y).\end{eqnarray} \end{prop}

\begin{rem} Formulas of \cite{FPR,FHR,Kun07,XW05} do not reconstruct $f(x)$ inside $S$ correctly if a part of support of $f$ lies outside $S$. This feature is absent in other reconstruction methods such as time reversal (the readers are referred to \cite{KKun,AKKun} for these discussion). I would be interesting to see whether formulas contained in \refcoro{fam} are any different in this regard.
\end{rem}

The paper is organized as follows. In \refsec{prgen}, we derive \reftheo{geninf} from \reftheo{main}. In \refsec{equ}, we show that \reftheo{genin} equivalent to \reftheo{geninf} and \reftheo{Kunt} is implied by \reftheo{geninf}. Two proofs of \reftheo{main} are presented in \refsec{mains}. In \refsec{spe}, Propositions \ref{fprop}, \ref{ffinite}, \ref{kuneq} and \ref{kunprop} are derived from \refcoro{fam} using the aforementioned choices of $\xi$ and $\varphi$. Proofs of some auxiliary results are given in \refsec{app}. Finally, some remarks are provided in \refsec{remm}.

\section{Derivation of \reftheo{geninf} from \reftheo{main}} \label{prgen}

Let $u$ solves (\ref{wave}) and $u_0$ be its even extension with respect to $t$. Then, $u_0$ solves the wave equation on $\RR^n \times \RR$. Therefore, \begin{equation} \label{fourier} \llg^2 \hat{u}_0(y,\llg) + \Delta \hat{u}_0(y,\llg) =0.\end{equation} Due to (\ref{he}) and (\ref{fourier}), we have the Green's identity: \begin{eqnarray}\label{grid} \hat{u}_0(x,\llg) = - \int\limits_{S} \left[\frac{\pdh \overline{G}(|x-y|,\llg)}{\pdh \nu_y} \hat{u}_0(y,\llg) - \frac{\pdh \hat{u}_0(y,\llg)}{\pdh \nu_y} \overline{G}(|x-y|,\llg)\right] d\sg(y).\end{eqnarray}
Since $f(x) =u_0(x,0)$, the Fourier inversion gives $$f(x) = \int\limits_\RR \hat{u}_0(x,\llg)d\llg.$$ Due to (\ref{grid}), we get \begin{eqnarray}\nonumber f(x) &=& - \int\limits_{S} \int\limits_{\RR} \frac{\pdh \overline{G}(|x-y|,\llg)}{\pdh \nu_y} \hat{u}_0(y,\llg) d\llg d\sg(y)\\  \label{f} &+& \int\limits_{S}\int\limits_{\RR} \frac{\pdh \hat{u}_0(y,\llg)}{\pdh \nu_y} \overline{G}(|x-y|,\llg) d\llg d\sg(y).\end{eqnarray}

Let $S(0,r)$ be the sphere centered at the origin with radius $r$. For any function $H \in C^1(\RR^n)$, by changing variables, we derive \begin{eqnarray*} \left. \frac{d}{dr}\left(r^{-n+1}\int\limits_{S(0,r)} H(y) d\sg(y)\right)\right|_{r=1} =\left. \frac{d}{dr}\left(\int\limits_{S} H(ry) d\sg(y)\right)\right|_{r=1}\end{eqnarray*}
That is, \begin{eqnarray*} \left. \frac{d}{dr}\left(r^{-n+1}\int\limits_{S(0,r)} H(y) d\sg(y)\right)\right|_{r=1} = \int\limits_{S} \frac{\pdh H}{\pdh \nu_y}(y) d\sg(y). \end{eqnarray*}
Applying this equality for $$H(y) = \int\limits_{\RR} \overline{G}(|y-x|,\llg) \hat{u}_0(y,\llg) d\llg,$$ we get \begin{eqnarray*} && \int\limits_{S} \frac{\pdh}{\pdh \nu_y}  \int\limits_{\RR} \overline{G}(|y-x|,\llg) \hat{u}_0(y,\llg) d\llg d\sg(y) = \\ &&\left. \frac{d}{dr}\left(r^{-n+1}\int\limits_{S(0,r)} \int\limits_{\RR} \overline{G}(|y-x|,\llg) \hat{u}_0(y,\llg) d\llg d\sg(y)\right)\right|_{r=1}. \end{eqnarray*}
 Due to Theorem \ref{main}, the right hand side is zero. Thus, \begin{eqnarray*} \int\limits_{S} \frac{\pdh}{\pdh \nu_y}  \int\limits_{\RR} \overline{G}(|y-x|,\llg) \hat{u}_0(y,\llg) d\llg d\sg(y) = 0. \end{eqnarray*}
Therefore, \begin{eqnarray}\label{mainap}\nonumber &&  \int\limits_{S}  \int\limits_{\RR} \frac{\pdh \overline{G}(|y-x|,\llg)}{\pdh \nu_y} \hat{u}_0(y,\llg) d\llg d\sg(y) \\ && = - \int\limits_{S}  \int\limits_{\RR} \overline{G}(|y-x|,\llg)\frac{ \pdh \hat{u}_0}{\pdh \nu_y}(y,\llg) d\llg d\sg(y). \end{eqnarray}
Combining this and (\ref{f}), we get
\begin{eqnarray}\label{gradin2}\nonumber f(x) &=& - 2 \int\limits_{S} \int\limits_{\RR} \frac{\pdh \overline{G}(|x-y|,\llg)}{\pdh \nu_y} \hat{u}_0(y,\llg) \\ &=& - 2 \int\limits_{S} \int\limits_{\RR} \overline{G}_s(|x-y|,\llg) \hat{u}_0(y,\llg) \frac{\left<y-x,y\right>}{|x-y|}.\end{eqnarray} Here $\overline{G}_s(s,\llg)$ is the derivative of $\overline{G}(s,\llg)$ with respect to $s$.

Applying Theorem \ref{main} for $\Og =B$, we get \begin{eqnarray} \label{range}\int\limits_{S} \int\limits_{\RR} \overline{G}(|x-y|,\llg) \hat{u}_0(y,\llg)=0.\end{eqnarray}
Taking the derivative with respect to $x$ of the above identity along the direction $\xi$, we obtain $$- 2\int\limits_{S} \int\limits_{\RR} \overline{G}_s(|x-y|,\llg) \hat{u}_0(y,\llg) \frac{\left<y-x,\xi \right>}{|x-y|} =0.$$
Subtracting this equality from (\ref{gradin2}), we conclude that \begin{eqnarray*} f(x) &=&- 2 \int\limits_{S} \int\limits_{\RR} \overline{G}_s(|x-y|,\llg)\hat{u}_0(y,\llg) \frac{\left<y-x,y-\xi \right>}{|x-y|}d\llg d\sg(y) \\ &=&- 2 \int\limits_{S}
\left.\left(\frac{d}{ds} \int\limits_{\RR} \overline{G}(s,\llg)\hat{g}_0(y,\llg) d\llg \right)\right|_{s=|x-y|}
\frac{\left<y-x,y-\xi \right>}{|x-y|} d\sg(y).\end{eqnarray*} Theorem \ref{geninf} is proved.\eproof

\section{Derivation of Theorems \ref{genin} and \ref{Kunt}}\label{equ}
We prove \reftheo{genin} by showing that it is equivalent to \reftheo{geninf}. Indeed, it suffices to prove the following result:
\begin{lemma} \label{finchlem} Suppose that $v \in C[0,\infty)$ such that $\hat{v}$ has proper decay at infinity, say $\hat{v}(\llg) = O(|\llg|^{-2})$ as $\llg \rightarrow \infty$. Then for any $s>0$, \begin{eqnarray}\label{Finchb} \int\limits_{\RR} \overline{G}(s,\llg)\hat{v}(\llg) d\llg = \mW(v)(s).\end{eqnarray}
\end{lemma}
In order to prove this lemma, we need the following explicit formula for $G$:
\begin{prop}\label{auxi}  Let $c_n$ be as in (\ref{cn}). For any $s>0$, we have \begin{eqnarray} \label{fund} G(s,\llg) = \left\{\begin{array}{l}
c_n \left(\frac{1}{s} \frac{d}{ds}\right)^{\frac{n-2}{2}} \left(\int\limits_{s}\limits^{\infty} \frac{e^{i \llg t}}{\sqrt{t^2-s^2}} dt \right),\mbox{ if $n$ is even},\\ c_n \left[\left(\frac{1}{s} \frac{d}{ds} \right)^{\frac{n-3}{2}}\frac{e^{i \llg s}}{s}\right],\mbox{ if $n$ is odd}. \end{array} \right.\end{eqnarray} \end{prop} This formula must be well known. For completeness, its proof is given in Section \ref{app}.

{\bf Proof of \reflemm{finchlem}} The relation \mref{fund} means that $G(s,\llg)=\mW(e^{i \llg t})(s)$. Keeping this fact in mind, one sees that the following proof consists of just changing the order of the operator $\mW$ with the integral sign:
\begin{itemize}
\item For even $n$, due to the decay of $\hat{v}$ and \refprop{auxi}, the following calculations are valid: \begin{eqnarray*} \int\limits_{\RR} \overline{G}(s,\llg)\hat{v}(\llg) d\llg &=& c_n \int\limits_{\RR} \left(\frac{1}{s} \frac{d}{ds} \right)^{\frac{n-2}{2}}\left(\int\limits_s\limits^\infty \frac{e^{-i \llg t}}{\sqrt{t^2-s^2}} dt \right) \hat{v}(\llg) d\llg\\ &=& c_n \left(\frac{1}{s} \frac{d}{ds} \right)^{\frac{n-2}{2}} \int\limits_s\limits^\infty \frac{1}{\sqrt{t^2-s^2}} \left(\int\limits_{\RR} e^{- i \llg t}\hat{v}(\llg) d \llg \right) dt\\&=& c_n \left(\frac{1}{s} \frac{d}{ds} \right)^{\frac{n-2}{2}} \int\limits_s\limits^\infty \frac{v(t)}{\sqrt{t^2-s^2}} dt = \mW(v)(s).\end{eqnarray*}

\item For an odd $n$, due to \refprop{auxi} and decay of $\hat{v}$, \begin{eqnarray*} \int\limits_{\RR} \overline{G}(s,\llg)\hat{v}(\llg) d\llg &=& c_n \int\limits_{\RR}\left(\frac{1}{s} \frac{d}{ds} \right)^{\frac{n-3}{2}}\left(\frac{e^{- is \llg }}{s} \right) \hat{v}(\llg) d\llg \\ &=& c_n\left(\frac{1}{s} \frac{d}{ds}\right)^{\frac{n-3}{2}}\left(\int\limits_{\RR}
\frac{e^{- i s \llg }}{s} \hat{v}(\llg) d\llg \right)\\ &=& c_n \left(\frac{1}{s} \frac{d}{ds} \right)^{\frac{n-3}{2}} \left(\frac{v(s)}{s}\right) = \mW(v)(s).\end{eqnarray*}
\end{itemize}
The lemma is proved.\eproof

We now prove the following result, which shows that \reftheo{geninf} implies \reftheo{Kunt}:
\begin{lemma}\label{kunlem}Let $f \in C_0^\infty(\RR^n)$, $g=\mT(f)$ and $g_0$ is its even extension with respect to $t$. Then \begin{eqnarray*}\int\limits_{\RR} \overline{G}(s,\llg)\hat{g}_0(y,\llg) d\llg =\frac{2}{\pi} \int\limits_0\limits^\infty \llg \mR(s,\llg) \int\limits_0\limits^\infty \mR_S(f)(y,r) \mI(r,\llg) dr.\end{eqnarray*} \end{lemma}
We need the following auxiliary result:
\begin{prop}\label{TR} Let $f \in C_0^\infty(\RR^n)$ and $g=\mT(f)$. Then \begin{eqnarray*} \int\limits_0\limits^\infty g(y,s) e^{i\llg s} ds = -i \llg\int\limits_0 \limits^\infty \mR_S(f)(y,s) G(s, \llg) dt.\end{eqnarray*} \end{prop}

\bproof For the expository purpose, we provide here two proofs of this proposition.
\begin{itemize}
\item[1)] Let $\bar{u}(x,\llg) = \int\limits_0\limits^\infty u(x,s)e^{i \llg s} ds$. Due to \mref{wave}, $\bar{u}(.,\llg)$ is the solution for \begin{eqnarray*} \Delta U(x) +\llg^2 U(x) = i \llg f(x),\end{eqnarray*} obtained by limiting absorption. Due to \mref{he}, we obtain \begin{eqnarray*} \bar{u}(y,\llg) &=& -i \llg \int\limits_{\RR^n} f(x) G(|y-x|,\llg) dx \\ &=& -i\llg \int\limits_0\limits^\infty (\mR_S f)(y,s) G(s,\llg) ds.\end{eqnarray*} This proves the proposition.
\item[2)] Consider the transform \begin{eqnarray*}\mB(v)(s) = \left\{\begin{array}{l} (-1)^{\frac{n-2}{2}}c_n \frac{d}{ds}\left(\frac{1}{s} \frac{d}{ds}\right)^{\frac{n-2}{2}} \left(\int\limits_0\limits^s \frac{v(t)}{\sqrt{s^2-t^2}} dt\right), \mbox{ if $n$ is even},\\ (-1)^{\frac{n-3}{2}}c_n \frac{d}{ds}\left(\frac{1}{s} \frac{d}{ds}\right)^{\frac{n-3}{2}} \left( \frac{v(s)}{s} \right), \mbox{ if $n$ is odd} \end{array} \right. \end{eqnarray*}
    Let $g \in C^\infty(S\times [0,\infty))$ , we also define $\mB(g)(y,s) = \mB(g_y)(s)$ where $g_y(s) =g(y,s)$. Let $g=\mT(f)$, we then have (e.g., \cite{CH2,EVB}) $g=\mB(\mR_S f).$ Therefore, \begin{eqnarray*} \int\limits_0\limits^\infty g(y,s) e^{i\llg s} ds = \int\limits_0\limits^\infty \mB(\mR_S f)(y,s) e^{i\llg s} ds.\end{eqnarray*} Substituting the expression of $\mB$ into the integral and integrating by parts, we obtain \begin{eqnarray*} \int\limits_0\limits^\infty g(y,s) e^{i\llg s} ds = -i \llg \int\limits_0\limits^\infty (\mR_S f)(y,s) \mW(e_\llg)(s) ds,\end{eqnarray*} where $e_\llg(s) =e^{i\llg s}$. Due to \refprop{auxi}, we have \begin{eqnarray*} \int\limits_0\limits^\infty g(y,s) e^{i\llg s} ds = -i \llg \int\limits_0\limits^\infty (\mR_S f)(y,s) G(s,\llg) ds.\end{eqnarray*} This completes the proof. \eproof
\end{itemize}

{\bf Proof of \reflemm{kunlem}} Due to \mref{fund}, $G(s,-\llg)=\overline{G}(s,\llg)$. Since $\hat{g}_0$ is even, we have \begin{eqnarray}\label{kl}\int\limits_{\RR} \overline{G}(s,\llg)\hat{g}_0(y,\llg) d\llg = 2 \int\limits_0\limits^\infty \mR(s,\llg)\hat{g}_0(y,\llg) d\llg.\end{eqnarray}
Since $g_0$ is the even extension of $g$ with respect to $t$, \begin{eqnarray*}\hat{g}_0(y,\llg) = \frac{1}{2 \pi} \int\limits_\RR g_0(y,t) e^{i\llg t} dt = \frac{1}{\pi} Re \int\limits_0\limits^\infty g(y,t) e^{i\llg t} dt.\end{eqnarray*} Due to \refprop{TR}, we get
\begin{eqnarray*}\hat{g}_0(y,\llg) &=& \frac{1}{\pi} Re \left( -i \llg\int\limits_0 \limits^\infty \mR_S(f)(y,s) G(s, \llg) dt \right) \\&=& \frac{\llg}{\pi} \int\limits_0 \limits^\infty \mR_S(f)(y,t) \mI(\llg,t) dt.\end{eqnarray*} From \mref{kl}, we arrive at \begin{eqnarray*}\int\limits_{\RR} \overline{G}(s,\llg)\hat{g}_0(y,\llg) d\llg =\frac{2}{\pi} \int\limits_0\limits^\infty \llg \mR(s,\llg) \int\limits_0\limits^\infty \mR_S(f)(y,r) \mI(r,\llg) dr.\end{eqnarray*} The proof is completed.
\eproof

\section{Proof of \reftheo{main}}\label{mains}
In this section, we present two proofs of Theorem \ref{main}. One of them relies on some results of \cite{FPR,FHR} and works in the time domain while the other one is self-contained and works in the frequency domain. Without loss of generality, we can assume that $\Og = B$, the open unit ball. We need to prove that if $f \in C^\infty(\overline{B})$ and $g=\mT(f)$ then for all $x \in B$ \begin{equation} \label{mainre} \int\limits_{S} \int\limits_{\RR} \overline{G}(|x-y|,\llg)\hat{g}_0(y,\llg) d\llg d\sg(y) =0,\end{equation} where $g_0$ is the even extension of $g$ with respect to $t$.

\subsection{The indirect proof} \label{indirect}
This proof is somewhat indirect since it uses inversion formulas in \cite{FPR,FHR}.

Due to \reflemm{finchlem}, equality (\ref{mainre}) is equivalent to $$\int\limits_{S} \mW(g_0)(y,|x-y|) d\sg(y) =0,$$ for all $x \in B$. Or, since $g=g_0|_{S \times [0,\infty)}$, \begin{equation} \label{mainin} \int\limits_{S} \mW(g)(y,|x-y|) d\sg(y) =0.\end{equation}

We introduce the following definition \begin{defi}Let $\widetilde{C}(S \times [0,\infty))$ be the space of all functions $h \in C^\infty(S \times [0,\infty))$ satisfying the following conditions: \begin{itemize} \item[i)] $h$ vanishes at $t=0$ to infinite order. \item[ii)] If $n$ is odd then $h$ is compactly supported. If $n$ is even then for any nonnegative integer $k$, $\|\pdh_t^k h(.,t)\|_{L^\infty(S)}=O(t^{-n-k+1})$ as $t \rightarrow \infty$.\end{itemize} \end{defi}

Then the operator $\mP$, introduced in \refsec{intro}, maps $C_0^\infty(\overline{B})$ into $\widetilde{C} (S \times [0,\infty))$. Indeed, one can show it by using the Kirchhoff-Poisson solution formulas for wave equation (e.g., \cite{EVB}).

We now prove the following auxiliary result:
\begin{lemma}\label{PW} Let $x \in B$ and $h \in \widetilde{C}(S \times [0,\infty))$. Then  $$\mP^*(h)(x) = \int\limits_{S}\mW(h)(y,|x-y|) d\sg(y),$$  where $\mP^*$ is the $L^2$-adjoint of $\mP$. \end{lemma} {\bf Proof} For an even $n$, a direct calculation (see \cite{FR07}) gives
\begin{eqnarray} \label{p*} \mP^*(h)(x) = c_n \int\limits_{S} \int\limits_{|x-y|}^\infty \frac{\left(\frac{d}{d t} \frac{1}{t} \right)^{\frac{n-2}{2}} h(y,t)}{\sqrt{t^2-|x-y|^2}} dt dy.\end{eqnarray} We observe that \begin{eqnarray*}\frac{1}{s} \frac{d}{ds} \int\limits_s\limits^\infty \frac{h(y,t)}{\sqrt{t^2-s^2}} dt &=& \frac{1}{s} \frac{d}{ds} \int\limits_s\limits^\infty \frac{h(y,t)}{t} \frac{t}{\sqrt{t^2-s^2}} dt\\ &=& - \frac{1}{s} \frac{d}{ds} \int\limits_s\limits^\infty \frac{d}{dt}
\left(\frac{h(y,t)}{t} \right) \sqrt{t^2-s^2}dt \\ & =& -\int\limits_s\limits^\infty \frac{d}{dt} \left(\frac{h(y,t)}{t} \right) \left(\frac{1}{s} \frac{d}{ds} \right)(\sqrt{t^2-s^2}) \\ &=& \int\limits_s\limits^\infty \frac{d}{dt} \left(\frac{h(y,t)}{t} \right)\frac{1}{\sqrt{t^2-s^2}}.\end{eqnarray*} Hence, by induction, for all $k \geq 0$, \begin{eqnarray*}\int\limits_s\limits^\infty \frac{\left(\frac{d}{dt} \frac{1}{t} \right)^k
h(y,t)} {\sqrt{t^2-s^2}} =\left(\frac{1}{s} \frac{d}{ds}\right)^{k} \int\limits_s\limits^\infty \frac{h(y,t)}{\sqrt{t^2-s^2}}.\end{eqnarray*} Therefore, due to (\ref{p*}), we have \begin{eqnarray*} \mP^*(h)(x) &=& c_n \int\limits_{S} \left(\frac{1}{s}
\frac{d }{d s} \right)^{\frac{n-2}{2}} \int\limits_{s}^\infty \frac{h(y,t)}{\sqrt{t^2-s^2}} dt dy \\ &=& \int\limits_{S} \mW(h)(y,|x-y|) d\sg(y).\end{eqnarray*}
For an odd $n$, a direct calculation (see \cite{FPR}) gives \begin{eqnarray*}\mP^*(g)(x) &=& c_n \int\limits_S\left. \left(\frac{1}{t}
\frac{d}{dt} \right)^{\frac{n-3}{2}}\left(\frac{h(y,t)}{t} \right)\right|_{t=|y-x|} d\sg(y) \\ &=& \int\limits_S \mW(h)(y,|x-y|)
d\sg(y).\end{eqnarray*} This proves the lemma.\eproof

This lemma tells us that (\ref{mainin}) is equivalent to $\mP^*(g) =0$. Equivalently,
\begin{equation}\label{sub} \mP^* \pdh_t \mP(f) =0.\end{equation}
We now derive this equality from inversion formula in \cite{FPR,FHR}. Indeed, if $n$ is even, \cite[Theorem 4]{FHR}\label{Fincheven} gives
\begin{eqnarray*} f(x) &=&  -2\left(\mP^*t\pdh_t^2 \mP f \right)(x), \\ f(x) &=& -2 \left(\mP^* \pdh_t t \pdh_t \mP f \right)(x).\end{eqnarray*}
By subtracting these two inversion formulas we obtain the desired equality (\ref{sub}).
If $n$ is odd, \cite[Theorem 3]{FPR}  gives \begin{eqnarray*}f(x) &=&  -2\left(\mP^* \pdh_t^2 t  \mP f \right)(x), \\ f(x) &=& -2 \left(\mP^* \pdh_t t \pdh_t \mP f \right)(x).\end{eqnarray*}
Again, subtracting these two inversion formulas, we obtain the desired equality (\ref{sub}). The proof of \reftheo{main} is completed.

\begin{rem}\label{rangerem} As shown in \cite{FR06}, identity (\ref{sub}) is the complete range description for the operator $\mP$ (or, equivalently, for $\mT$) when $n$ is odd.
 \end{rem}

\subsection{The direct proof}\label{direct}
Let $u_1$ be the extension of $u$ by zero for $t<0$, $\hat{u}_1$ be its Fourier transform, and $g_1= u_1|_{S \times \RR}$. Applying Lemma \ref{finchlem}, we have $$\int\limits_{S} \int\limits_{\RR} \overline{G}(|x-y|,\llg)\hat{g}_0(y,\llg)d \llg d\sg(y) = \int\limits_{S} \mW(g_0)(y,|x-y|) d\sg(y),$$ and $$\int\limits_{S} \int\limits_{\RR} \overline{G}(|x-y|,\llg)\hat{g}_1(y,\llg) d\llg d\sg(y)=\int\limits_{S} \mW(g_1)(y,|x-y|) d\sg(y).$$
Since $g_0=g_1$ on $S \times [0,\infty)$, one has $\mW(g_0)(y,|x-y|)=\mW(g_1)(y,|x-y|)$. The above two equalities then give \begin{eqnarray}\label{u0u1}\int\limits_{S} \int\limits_{\RR} \overline{G}(|x-y|,\llg)\hat{g}_0(y,\llg)d \llg d\sg(y) = \int\limits_{S} \int\limits_{\RR} \overline{G}(|x-y|,\llg)\hat{g}_1(y,\llg) d\llg d\sg(y). \end{eqnarray}
Due to \reflemm{kunlem}, \begin{eqnarray}\label{fouru1} \nonumber \hat{g}_1(y,\llg) &&= \frac{1}{2 \pi} \int\limits_0\limits^\infty g(y,s) e^{i \llg s} ds = \frac{-i \llg}{2 \pi} \int\limits_0^\infty \mR_S(f)(y,s) G(s,\llg)ds \\ &&= \frac{-i \llg}{2 \pi} \int\limits_{\RR^n} f(z) G(|z-y|,\llg)dz. \end{eqnarray} From (\ref{u0u1}) and (\ref{fouru1}), we get \begin{eqnarray*}&& \int\limits_{S} \int\limits_{\RR} \overline{G}(|x-y|,\llg)\hat{g}_0(y,\llg) d\llg d\sg(y)  \\ &=& \frac{-i}{2 \pi}\int\limits_{S} \int\limits_{\RR} \overline{G}(|x-y|,\llg) \llg \int\limits_{\RR^n}G(|y-z|,\llg)f(z) dz d\llg d\sg(y) \\ &=&\frac{-i}{2 \pi} \int_{B} f(z) \left(\int\limits_{\RR} \llg \int\limits_{S} \overline{G}(|x-y|,\llg)G(|y-z|,\llg)d\sg(y) d\llg \right) d z. \end{eqnarray*} That is, \begin{eqnarray}\label{h} \int\limits_{S} \int\limits_{\RR} \overline{G}(|x-y|,\llg)\hat{g}_0(y,\llg) d\llg d\sg(y) =\frac{- i}{2 \pi} \int\limits_{B} f(z) \left(\int\limits_{\RR} \llg \varphi(\llg) d\llg \right) d z, \end{eqnarray}  where $$\varphi(\llg) = \int\limits_{S} \overline{G}(|x-y|,\llg) G(|y-z|,\llg) d\sg(y).$$
From the explicit formula for $G$ in Proposition \ref{auxi}, we see that $G(s,-\llg) = \overline{G}(s,\llg)$. Thus, \begin{equation} \label{con}\varphi(-\llg)= \overline{\varphi}(\llg).\end{equation}  Let $\mR(s,\llg) =Re(G(s,\llg)), \mI(s,\llg) =Im (G(s,\llg))$, and $$K(x,z,\llg) = \int\limits_{S} \mR(|x-y|,\llg) \mI(|z-y|, \llg) d\sg(y).$$ We claim the following symmetry whose proof is presented in Section \ref{app}: \begin{equation}\label{sym} K(x,z,\llg) = K(z,x,\llg) \mbox{ for any $x,z \in B$ and $\llg \in \RR$.}\end{equation} Assuming this symmetry, one has $Im (\varphi(\llg)) = K(x,z,\llg) - K(z,x,\llg)=0$. From (\ref{con}), we get $\varphi(-\llg)= \varphi(\llg)$. This implies the following equality: $$\int\limits_{\RR} \llg \varphi(\llg) d\llg =0.$$ Due to (\ref{h}), one concludes that $$\int\limits_{S} \int\limits_{\RR} \overline{G}(|x-y|,\llg)\hat{g}_0(y,\llg) d\llg d\sg(y) =0.$$ The proof is completed.

\section{Some special cases} \label{spe}

In this section, we derive Propositions \ref{fprop}, \ref{ffinite}, \ref{kuneq} and \ref{kunprop} from \refcoro{fam} by proper choices of $\xi$ and $\varphi$.
\subsection{Propositions \ref{fprop}}
By choosing $\xi=x$ and $\varphi(x) =c$ in \mref{fam2}, we obtain \begin{eqnarray*}f(x)  =- 2 \int\limits_{S} \left.\pdh_s\mW(g)(y,s)\right|_{s=|x-y|}|x-y| d\sg(y) + c \int\limits_{S} \mW(g)(y,|x-y|)d\sg(y).\end{eqnarray*} Here we use the notation $\pdh_s$ for $\frac{d}{ds}$.
That is, \begin{eqnarray}\label{Ftype}f(x)=- 2 \int\limits_{S} \left[s \pdh_s \mW g\right](y,|x-y|)d\sg(y)+c \int\limits_{S} \mW(g)(y,|x-y|)d\sg(y).\end{eqnarray}
We now claim an identity, whose proof can be found in Section \ref{app}: \begin{equation}\label{simp} s \pdh_s \mW (g) =  \mW (s \pdh_s g) - (n-2) \mW(g).\end{equation}
Assuming this claim, we see that (\ref{Ftype}) is equivalent to \begin{eqnarray*} f(x) &=& -2 \int\limits_{S}\mW (s \pdh_s g)(y,|x-y|) d\sg(y) \\ &+& [2(n-2)+c]\int\limits_{S}\mW (g)(y,|x-y|) d\sg(y).\end{eqnarray*}
Due to \reflemm{PW}, we obtain \begin{eqnarray}\label{fr} f(x) = -2 \mP^* \left(s \pdh_s g\right)(x) + [2(n-2)+c]\mP^*(g)(x)
.\end{eqnarray}
\begin{itemize}
\item If $n=-2(n-2)$ then \mref{fr} becomes \begin{eqnarray*} f(x) = -2 \mP^* \left(s \pdh_s g\right)(x) = -2 (\mP^*s\pdh_s^2 \mP f)(x) \end{eqnarray*}
\item If $n=-2(n-1)$ then \mref{fr} becomes \begin{eqnarray*} f(x) = -2 \mP^* \left(s \pdh_s g\right)(x) - 2\mP^*(g)(x) = -2 (\mP^* \pdh_s s \pdh_s \mP f)(x) \end{eqnarray*}
\item If $n=-2n$ then \mref{fr} becomes \begin{eqnarray*} f(x) = -2 \mP^* \left(s \pdh_s g\right)(x) - 4\mP^*(g)(x)= -2 (\mP^*\pdh_s^2 s\mP f)(x)
.\end{eqnarray*}
\end{itemize}
Propositions \ref{fprop} is proved (in this proof we have used the variable $s$ in place of $t$).

\subsection{Proposition \ref{ffinite}}

Choosing $\varphi(x) =0$ in \mref{fam1}, we obtain
\begin{eqnarray}\label{finite} f(x) =- \frac{4}{\pi} \int\limits_{S} \left(\frac{d }{ds}K_n(y,s) \right)_{s=|x-y|} \frac{\left<y-x,y- \xi\right>}{|x-y|} d\sg(y),\end{eqnarray}
Due to (\ref{green}) we have $G(s,\llg) = \frac{\llg^{n-2}}{4 (2 \pi)^{\frac{n-2}{2}}} \left[i J(\llg s) - N(\llg s) \right]$. Thus, \begin{eqnarray}\label{connect} \mR(s,\llg) = -\frac{\llg^{n-2}}{4 (2 \pi)^{\frac{n-2}{2}}} N(\llg s), ~ \mI(r,\llg) = \frac{\llg^{n-2}}{4 (2 \pi)^{\frac{n-2}{2}}} J(\llg s) .\end{eqnarray}
Hence, \begin{eqnarray*} K_n(y,s) &=&  \int\limits_0\limits^\infty \llg \mR(s\llg) \int\limits_0\limits^\infty \mR_S(f)(y,r) \mI(r \llg) dr d\llg\\ &=& \frac{-1}{16 (2\pi)^{n-2}}\int\limits_0\limits^\infty \llg^{2n-3} N(s\llg) \int\limits_0\limits^\infty \mR_S(f)(y,r) J(r \llg) dr d\llg.\end{eqnarray*}
Since $f$ is supported inside $\overline{B}$, we have $\mR_S(f)(y,r)=0$ for $r \geq 2$. Therefore, \begin{eqnarray*} K_n(y,s) &=&  \frac{-1}{16 (2\pi)^{n-2}}\int\limits_0\limits^\infty \llg^{2n-3} N(s\llg) \int\limits_0\limits^2 \mR_S(f)(y,r) J(r \llg) dr d\llg \\ &=& \frac{-1}{16 (2\pi)^{n-2}}k_n(y,s),\end{eqnarray*} where, as defined in \refsec{intro}, \begin{eqnarray*}k_n(y,s) = \int\limits_0\limits^\infty \llg^{2n-3} N(s\llg) \int\limits_0\limits^2 \mR_S(y,r) J(r \llg) dr d\llg.\end{eqnarray*}

From \mref{finite}, we arrive at
\begin{eqnarray}\label{last} f(x) =\frac{1}{2 (2\pi)^{n-1}} \int\limits_{S} \left(\frac{d }{ds}k_n(y,s) \right)_{s=|x-y|} \frac{\left<y-x,y- \xi\right>}{|x-y|} d\sg(y),\end{eqnarray}
Choosing $\xi=x$, we have \begin{eqnarray*} f(x) &=& \frac{1}{2 (2\pi)^{n-1}} \int\limits_{S} \left(\frac{d}{ds}k_n(y,s) \right)_{s=|x-y|} |x-y| d\sg(y).\end{eqnarray*}
That is, \begin{eqnarray}\label{fe} f(x) &=& \frac{1}{2 (2\pi)^{n-1}} \int\limits_{S} \left(s \frac{d}{ds}k_n(y,s) \right)_{s=|x-y|} d\sg(y).\end{eqnarray}
 We claim that for even $n$ \begin{eqnarray} \label{evenf} k_n(y,s) = \frac{(-1)^{\frac{n-2}{2}}}{\pi}  \int\limits_0\limits^2 \left(\frac{d}{dr} \frac{1}{r}\right)^{n-1}\mR_S(f)(y,r) \ln |r^2-s^2| dr.\end{eqnarray} Here the integral is understood in the principal value sense. A proof of this claim will be given in \refsec{apevenf}. Assuming it, we obtain \begin{eqnarray*}s \frac{d}{ds}k_n(y,s) &=& \frac{2(-1)^{\frac{n-2}{2}}}{\pi} \int\limits_0\limits^2 \left(\frac{d}{dr} \frac{1}{r}\right)^{n-1}\mR_S(f)(y,r) \frac{s^2}{s^2-r^2} dr \\ &=& \frac{2(-1)^{\frac{n-2}{2}}}{\pi} \int\limits_0\limits^2 \left(\frac{d}{dr} \frac{1}{r}\right)^{n-1}\mR_S(f)(y,r) \left[1+\frac{r^2}{s^2-r^2}\right] dr.\end{eqnarray*}
Since \begin{eqnarray*}\int\limits_0\limits^2 \left(\frac{d}{dr} \frac{1}{r}\right)^{n-1}\mR_S(f)(y,r) dr = \left.\frac{1}{r}\left(\frac{d}{dr} \frac{1}{r}\right)^{n-2}\mR_S(f)(y,r)\right|_0^2 =0,\end{eqnarray*} we get \begin{eqnarray*}s \frac{d}{ds}k_n(y,s) &=& \frac{2(-1)^{\frac{n-2}{2}}}{\pi} \int\limits_0\limits^2 \left(\frac{d}{dr} \frac{1}{r}\right)^{n-1}\mR_S(f)(y,r) \frac{r^2}{s^2-r^2}dr \\ &=&  \frac{2(-1)^{\frac{n}{2}}}{\pi} \int\limits_0\limits^2 \left(\frac{d}{dr} \frac{1}{r}\right)^{n-1}\mR_S(f)(y,r) \frac{r^2}{r^2-s^2}dr.\end{eqnarray*}
Integrating by part (which can be justified for the principal value integral in question), we obtain \begin{eqnarray*}s \frac{d}{ds}k_n(y,s)  &=&  \frac{(-1)^{\frac{n-2}{2}}}{\pi} \int\limits_0\limits^2 \frac{d}{dr}r\left(\frac{d}{dr} \frac{1}{r}\right)^{n-1}\mR_S(f)(y,r) \ln|r^2-s^2|dr.\end{eqnarray*}
From \mref{fe}, we arrive at \begin{eqnarray*}f(x) = \frac{(-1)^{\frac{n-2}{2}}}{(2\pi)^{n}} \int\limits_{S} \int\limits_0\limits^2 \left[\frac{d}{dr}r\left(\frac{d}{dr} \frac{1}{r}\right)^{n-1}\mR_S(f)\right](y,r) \ln \left|r^2-|x-y|^2\right| dr d\sg(y).\end{eqnarray*}
Equivalently, $f(x)$ is equal to \begin{eqnarray*} \frac{(-1)^{\frac{n-2}{2}}\og_n}{(2\pi)^{n}} \int\limits_{S} \int\limits_0\limits^2 \left[\frac{d}{dr}r\left(\frac{d}{dr} \frac{1}{r}\right)^{n-1} r^{n-1} \mM_S(f)\right](y,r) \ln \left|r^2-|x-y|^2\right| dr d\sg(y).\end{eqnarray*}
\eproof

\subsection{Propositions \ref{kuneq} and \ref{kunprop}} \label{Kunsub}

Choosing $\xi=0$ in (\ref{last}), we obtain \begin{eqnarray*} f(x) = \frac{1}{2 (2\pi)^{n-1}} \int\limits_{S} \left(\frac{d }{ds}k_n(y,s) \right)_{s=|x-y|} \frac{\left<y-x,y\right>}{|x-y|} d\sg(y),\end{eqnarray*}
That is, \begin{eqnarray}\label{mine} f(x) = \frac{- 1}{2 (2\pi)^{n-1}} \nabla_x \int\limits_{S} n(y) k_n(y,|x-y|) d\sg(y).\end{eqnarray}
This proves \refprop{kuneq}. In order to prove \refprop{kunprop}, we need only show that $h(y,s) = 2k_n(y,s)$ for any $s>0$. Indeed, it is a consequence of the following lemma whose proof can be found in Section \ref{app}:
\begin{lemma}\label{toprk} Suppose that $h \in C^\infty[0,\infty)$ such that $h$ does not grow too fast at infinity and $h^{(i)}(0)=0$ for all $0 \leq i \leq n-3$ if $n$ is odd, and $h^{(i)}(0)=0$ for all $0 \leq i \leq n-2$ if $n$ is even. Then \begin{eqnarray}\label{eqkun} \nonumber \int\limits_0\limits^\infty \llg^{2n-3} N(s\llg) \int\limits_0\limits^\infty h(r) J(r \llg) dr d\llg = - \int\limits_0\limits^\infty \llg^{2n-3} J(s\llg) \int\limits_0\limits^\infty h(r) N(r \llg)d\llg.\end{eqnarray} \end{lemma}

\section{Proofs of auxiliary statements}\label{app}
\subsection{Proof of Proposition \ref{auxi}}
We now prove \begin{eqnarray}\label{greenr} G(s,\llg) = \left\{\begin{array}{l}
c_n \left(\frac{1}{s} \frac{d}{ds}\right)^{\frac{n-2}{2}} \left(\int\limits_{s}\limits^{\infty} \frac{e^{i \llg t}}{\sqrt{t^2-s^2}} dt \right),\mbox{ if $n$ is even},\\ c_n \left(\frac{1}{s} \frac{d}{ds} \right)^{\frac{n-3}{2}}\left(\frac{e^{i \llg s}}{s}\right),\mbox{ if $n$ is odd}. \end{array} \right.\end{eqnarray}

Indeed, from (\ref{green}), we get \begin{equation} \label{ag} G(s,\llg)= \frac{i}{4} \left(\frac{\llg}{2 \pi s}
\right)^{\frac{n-2}{2}} H_{\frac{n-2}{2}}^{(1)}(\llg s) = \frac{i}{4} \frac{\llg^{n-2}}{(2 \pi)^{\frac{n-2}{2}}}
\frac{H^{(1)}_{\frac{n-2}{2}}(\llg s)}{(\llg s)^{\frac{n-2}{2}}},\end{equation} where $H^{(1)}_{\frac{n-2}{2}}$ is the Hankel function of the first kind.
\begin{itemize}

%\item Consider $n=2$. It is known that (e.g. \cite [page 180]{Watson}) $$H_0^{(1)}(s)= \frac{-2i}{\pi} \int_1^\infty \frac{e^{ist}}{\sqrt{t^2-1}} dt.$$ From (\ref{ag}), we have $$G(s,\llg) = \frac{i}{4} H_0^{(1)}(\llg s) =
%\frac{1}{2\pi} \int_1^\infty \frac{e^{i \llg st}}{\sqrt{t^2-1}} dt.$$ This confirms (\ref{greenr}) for $n=2$.

\item Let $n$ be an even number.  Using the equality (e.g. \cite[page 74]{Watson}) \begin{equation}\label{wat} \frac{H^{(1)}_{\nu+m}(s)}{s^{\nu+m}}=(-1)^m \left(\frac{1}{s} \frac{d}{ds}\right)^{m} \left(\frac{H_\nu^{(1)}(s)}{s^{\nu}}\right),\end{equation} we obtain \begin{eqnarray*}\frac{H^{(1)}_{\frac{n-2}{2}}(\llg s)}{(\llg s)^{\frac{n-2}{2}}} &=& (-1)^{\frac{n-2}{2}}\left. \left[\left(\frac{1}{s} \frac{d}{ds}\right)^{\frac{n-2}{2}} H_0^{(1)}(s)\right]\right|_{\llg s} \\ &=& (-1)^{\frac{n-2}{2}} \llg^{-(n-2)} \left(\frac{1}{s} \frac{d}{ds}\right)^{\frac{n-2}{2}} H_0^{(1)}(\llg s). \end{eqnarray*}
Hence, (\ref{ag}) gives
 \begin{eqnarray*}G(s,\llg) &=& \frac{i}{4} \frac{(-1)^{\frac{n-2}{2}}}{(2 \pi)^{\frac{n-2}{2}}}\left(\frac{1}{s} \frac{d}{ds}\right)^{\frac{n-2}{2}}
H_0^{(1)}(\llg s).\end{eqnarray*} Since (e.g., \cite [page 170]{Watson}) $$H_0^{(1)}(s)= \frac{-2i}{\pi} \int_1^\infty \frac{e^{ist}}{\sqrt{t^2-1}} dt= \frac{-2i}{\pi} \int_s^\infty \frac{e^{it}}{\sqrt{t^2-s^2}} dt,$$ we conclude that $$G(s,\llg) = \frac{(-1)^{\frac{n-2}{2}}}{(2 \pi)^{\frac{n}{2}}} \left(\frac{1}{s}
\frac{d}{ds}\right)^{\frac{n-2}{2}} \left(\int_s^\infty \frac{e^{i \llg t}}{\sqrt{t^2-s^2}} dt \right).$$ This confirms (\ref{greenr}) for even $n$.

%\item Now, consider the case $n=3$. It is known that  (e.g. \cite[page 170]{Watson}) $$H_{1/2}^{(1)}(s) = - i \sqrt{\frac{2}{\pi s}} e^{is}.$$ Thus, (\ref{ag}) gives \begin{eqnarray*}G(s,\llg) = \frac{i}{4} \left( \frac{\llg}{2 \pi s} \right)^{1/2} H^{(1)}_{1/2}(\llg s) = \frac{1}{4\pi} \frac{e^{i \llg s}}{s}. \end{eqnarray*}

\item Let $n$ be an odd number. Using (\ref{wat}) again, one gets \begin{eqnarray*}\frac{H^{(1)}_{\frac{n-2}{2}}(\llg s)}{(\llg s)^{\frac{n-2}{2}}}=(-1)^{\frac{n-3}{2}} \llg^{-(n-3)} \left[\left(\frac{1}{s} \frac{d}{ds} \right)^{\frac{n-3}{2}}\frac{H^{(1)}_{\frac{1}{2}}(\llg s)}{(\llg s)^{\frac{1}{2}}}\right].\end{eqnarray*}

Hence, from (\ref{ag}), we obtain \begin{eqnarray*}G(s,\llg) &=& \frac{i}{4} \frac{(-1)^{\frac{n-3}{2}}\llg}{(2 \pi) ^{\frac{n-2}{2}}} \left[\left(\frac{1}{s} \frac{d}{ds}
\right)^{\frac{n-3}{2}}\frac{H^{(1)}_{\frac{1}{2}}(\llg s)}{(\llg s)^{\frac{1}{2}}}\right].\end{eqnarray*}

 Since (e.g., \cite[page 487]{CH1}) $$H_{\frac{1}{2}}^{(1)}(s) = - i \sqrt{\frac{2}{\pi s}} e^{is},$$ we conclude that $$G(s,\llg) = \frac{(-1)^{\frac{n-3}{2}}}{2(2 \pi) ^{\frac{n-1}{2}}} \left[\left(\frac{1}{s} \frac{d}{ds} \right)^{\frac{n-3}{2}}\frac{e^{i \llg s}}{s}\right].$$  This confirms (\ref{greenr}) for odd $n$.
 \end{itemize}

\subsection{Proof of identity (\ref{sym})}

We prove that for all $\llg \in \RR$, the function $$K(x,z,\llg) = \int\limits_{S} \mR(|x-y|,\llg) \mI(|z-y|, \llg) d\sg(y),$$ is symmetric with respect to $x,z \in B$. Here we follow the line of reasoning used in \cite{Kun07}.

Due to (\ref{connect}), we get \begin{equation} \label{connect2} K(x,z,\llg) = c_0 \int\limits_{S} N(\llg |x-y|) J(\llg|z-y|) d\sg(y),\end{equation} where $c_0$ is a constant depending only on $n$ and $\llg$. We recall the following identities from \cite{Kun07} for $|x|>r_0$: \begin{eqnarray} \label{Kunid1} \int_{|y|=r_0} Y_l^{k}(\hat{y})J(\llg |x-y|)dy &=& c_1 J_{k}(\llg r_0) J_{k}(\llg|x|) Y_l^{k}(\hat{x}),
\\ \label{Kunid2} \int_{|y|=r_0} Y_l^{k}(\hat{y}) N(\llg|x-y|) dy &=& c_1 J_{k}(\llg r_0) N_{k}(\llg|x|) Y_l^{k}(\hat{x}), \end{eqnarray} where $\hat{x}=\frac{x}{|x|}$ and $Y_l^k$ is a spherical harmonic of order $k$. Consider the spherical harmonic expansion $$K(x,z,\llg)=  \sum_{(k,l),(k',l')} a^{k,l}_{k',l'}(\ag,\bg) Y_l^{k}(\hat{x})Y_{l'}^{k'}(\hat{z}),$$ where $\ag=|x|$ and $\bg=|z|$.

Due to (\ref{connect2}), we obtain \begin{eqnarray*} &&a^{k,l}_{k',l'}(\ag,\bg)\\ &&=\int\limits_S \int\limits_S Y_l^k(\hat{x}) Y_{l'}^{k'}(\hat{z}) K(\ag \hat{x}, \bg \hat{z}, \llg) d\sg(\hat{x}) d \sg(\hat{z}) \\ &&= c_0 \int\limits_S \int\limits_S Y_l^k(\hat{x}) Y_{l'}^{k'}(\hat{z}) \int\limits_S N(\llg |y-\ag \hat{x}|) J(\llg |y- \bg \hat{z}|) d \sg(y) d \hat{x} d \hat{z} \\ &&= c_0 \int\limits_S \left(\int\limits_S Y_l^k(\hat{x})N(\llg |y-\ag \hat{x}|)d \hat{x} \right) \left( \int\limits_S Y_{l'}^{k'}(\hat{z})J(\llg |y- \bg \hat{z}|)  d \hat{z}\right) d \sg(y).\end{eqnarray*}
Applying (\ref{Kunid1}) and (\ref{Kunid2}), we arrive at \begin{eqnarray*}a^{k,l}_{k',l'}(\ag,\bg)  &=& c_0c_1^2 J_k(\llg \ag) J_{k'}(\llg \bg) J_{k}(\llg) N_{k'}(\llg) \int\limits_S Y_k^l(y) Y_{k'}^{l'}(y) d\sg(y) \\ &=& c_0c_1^2 J_k(\llg \ag) J_{k'}(\llg \bg) J_{k}(\llg) N_{k'}(\llg) \delta_{k,k'} \delta_{l,l'}.\end{eqnarray*} Hence, $a^{k,l}_{k,l}(\ag,\bg) =a^{k,l}_{k,l}(\bg,\ag)$ and $K(x,z,\llg)=  \sum_{(k,l)} a^{k,l}_{k,l}(\ag,\bg) Y_l^{k}(\hat{x})Y_{l}^{k}(\hat{z})$. These two equalities give $K(x,z,\llg) =K(z,x,\llg)$. The proof is completed.

\subsection{Proof of identity (\ref{simp})}
We now prove the identity \begin{equation} \label{simpr} \left[s \frac{d}{ds} \mW \right](g) = \left[\mW s \frac{d}{d s}\right](g) - (n-2) \mW(g).\end{equation}
We first recall from (\ref{mw}) \begin{eqnarray}\label{mwr} \mW(g)(y,s) :=\left\{\begin{array}{l} c_n \left(\frac{1}{s} \frac{d}{ds} \right)^{\frac{n-2}{2}} \int\limits_{s}^\infty \frac{g(y,t)}{\sqrt{t^2-s^2}}  dt, \mbox{ if $n$ is even}, \\ c_n \left(\frac{1}{s} \frac{d}{ds} \right)^{\frac{n-3}{2}} \left(\frac{g(y,s)}{s}\right),\mbox{ if $n$ is odd}. \end{array} \right.\end{eqnarray}

\begin{itemize}
\item For $n=2$, due to (\ref{mwr}), we have  \begin{eqnarray*} \left(s \frac{d}{d s}\mW\right)(g)(y,s)= c_2 \left(s \frac{d}{d s}\right) \int\limits_s\limits^\infty \frac{g(y,t)}{\sqrt{t^2-s^2}} dt.\end{eqnarray*}

Taking integration by parts, one gets \begin{eqnarray*} \int\limits_s\limits^\infty \frac{g(y,t)}{\sqrt{t^2-s^2}} dt &=& - \int\limits_s\limits^\infty \frac{d}{dt} \left(\frac{g(y,t)}{t}\right)\sqrt{t^2-s^2} dt.\end{eqnarray*}

Hence, \begin{eqnarray*} && s \frac{d}{ds} \int\limits_s\limits^\infty
\frac{g(y,t)}{\sqrt{t^2-s^2}} dt = - s \int\limits_s\limits^\infty \frac{d}{dt} \left(\frac{g(y,t)}{t}\right)
 \frac{d}{ds}\left(\sqrt{t^2-s^2}\right) dt \\ &=& \int\limits_s\limits^\infty \frac{d}{dt} \left(\frac{g(y,t)}{t}\right) \frac{s^2}{\sqrt{t^2-s^2}} dt \\ &=& \int\limits_s\limits^\infty \frac{d}{dt} \left(\frac{g(y,t)}{t}\right)
\frac{t^2}{\sqrt{t^2-s^2}} dt - \int\limits_s\limits^\infty \frac{d}{dt} \left(\frac{g(y,t)}{t}\right) \sqrt{t^2-s^2} dt \end{eqnarray*}

Taking integration by parts again, we arrive at \begin{eqnarray*} s \frac{d}{ds} \int\limits_s\limits^\infty
\frac{g(y,t)}{\sqrt{t^2-s^2}} dt &=& \int\limits_s\limits^\infty \frac{d}{dt} \left(\frac{g(y,t)}{t}\right) \frac{t^2}{\sqrt{t^2-s^2}} dt \\ &+& \int\limits_s\limits^\infty \left(\frac{g(y,t)}{t}\right) \frac{t} {\sqrt{t^2-s^2}} dt.\end{eqnarray*}

Simplifying the right hand side, we get \begin{eqnarray}\label{ex2} s \frac{d}{ds} \int\limits_s\limits^\infty \frac{g(y,t)}{\sqrt{t^2-s^2}} dt = \int\limits_s\limits^\infty \frac{tg_t(y,t)}{\sqrt{t^2-s^2}} dt.\end{eqnarray}

Therefore, \begin{eqnarray*}\left(s \frac{d}{d s}\mW\right)(g)(y,s)&=& c_2 \int\limits_s\limits^\infty \frac{tg_t(y,t)}{\sqrt{t^2-s^2}} dt = \left(\mW s\frac{d}{ds}\right)g(y,s).\end{eqnarray*} This confirms (\ref{simpr}) for $n=2$.

\item Let $n>2$ be even. We first observe the simple identity $$\left(s \frac{d}{ds} \right) \left( \frac{1}{s}
\frac{d}{ds}\right) = \left(\frac{1}{s} \frac{d}{ds}\right)\left(s\frac{d}{ds} \right) -2 \left( \frac{1}{s} \frac{d}{ds} \right).$$

By induction, we get \begin{equation} \label{inter}\left(s \frac{d}{ds} \right) \left(\frac{1}{s} \frac{d}{ds}\right)^k = \left(\frac{1}{s}
\frac{d}{ds}\right)^k \left(s \frac{d}{ds} \right) -2k \left(\frac{1}{s} \frac{d}{ds} \right)^k.\end{equation}

Therefore, due to the definition of $\mW$ in (\ref{mwr}), \begin{eqnarray*} \left(s\frac{d}{d s}\mW\right)(g) &=& c_n \left(s \frac{d}{d s}\right)
\left(\frac{1}{d} \frac{d}{ds}\right)^{\frac{n-2}{2}} \int\limits_s\limits^\infty \frac{g(y,t)}{\sqrt{t^2-s^2}} dt \\ &=& c_n \left(\frac{1}{d}
\frac{d}{ds}\right)^{\frac{n-2}{2}} \left(s \frac{d}{d s}\right)\int\limits_s\limits^\infty \frac{g(y,t)}{\sqrt{t^2-s^2}} dt \\ &-& (n-2) c_n \left(\frac{1}{d} \frac{d}{ds}\right)^{\frac{n-2}{2}} \int\limits_s\limits^\infty \frac{g(y,t)}{\sqrt{t^2-s^2}} dt.\end{eqnarray*}
Due to (\ref{ex2}), we obtain \begin{eqnarray*} \left(s \frac{d}{d s}\right)\mW(g)(y,s)&=& c_n \left(\frac{1}{d} \frac{d}{ds}\right)^{\frac{n-2}{2}} \int\limits_s\limits^\infty \frac{tg_t(y,t)}{\sqrt{t^2-s^2}} dt \\ &-&(n-2) c_n \left(\frac{1}{d}\frac{d}{ds}\right)^{\frac{n-2}{2}}  \int\limits_s\limits^\infty \frac{g(y,t)}{\sqrt{t^2-s^2}} dt \\ &=& \left(\mW s\frac{d}{ds}\right)(g)(y,s) -(n-2) \mW(g)(y,s).\end{eqnarray*} This confirms (\ref{simpr}) for any even $n$.

%\item{Let n=3.} From the definition of $\mW$ in (\ref{mwr}), one derives \begin{eqnarray*} \left(s \frac{d}{ds}\mW\right)(g)(y,s) &=& c_3 \left(s \frac{\pdh}{\pdh s}\right) \frac{g(y,s)}{s}=  c_3 \left(u_s(y,s) - \frac{g(y,s)}{s} \right) \\ &=& \mW(sg_s) - \mW(g).\end{eqnarray*}

\item Let $n$ be odd. Due to (\ref{mwr}), \begin{eqnarray*} \left(s \frac{d}{ds}\mW\right)(g)(y,s)= c_n \left(s \frac{d}{ds}\right) \left(\frac{1}{s} \frac{d}{ds}\right)^{\frac{n-3}{2}} \left(\frac{g(y,s)}{s} \right).\end{eqnarray*}

This and (\ref{inter}) give \begin{eqnarray*}\left(s \frac{d}{ds}\mW\right)(g)(y,s)&=& c_n \left(\frac{1}{s} \frac{d}{ds}\right)^{\frac{n-3}{2}} \left(s \frac{d}{d s}\right) \left(\frac{g(y,s)}{s} \right) \\ &-& c_n (n-3)\left(\frac{1}{s} \frac{d}{ds}\right)^{\frac{n-3}{2}} \left(\frac{g(y,s)}{s} \right).\end{eqnarray*}

Using the identity $\left(s \frac{d}{d s}\right) \left(\frac{g(y,s)}{s} \right)= g_s(y,s) - \frac{g(y,s)}{s}$ for the first term of the right hand side of the above equality, we obtain
\begin{eqnarray*} \left(s \frac{d}{ds}\mW\right)(g)(y,s) &=& c_n \left(\frac{1}{s} \frac{d}{ds}\right)^{\frac{n-3}{2}} g_s(y,s)  \\ &-& c_n (n-2)\left(\frac{1}{s} \frac{d}{ds}\right)^{\frac{n-3}{2}}\left(\frac{g(y,s)}{s} \right) \\ &=& \left(\mW s \frac{d}{ds}\right)(g)(y,s) - (n-2)\mW(g)(y,s) .\end{eqnarray*} This finishes the proof of identity (\ref{simpr}).
\end{itemize}

\subsection{Proof of \mref{evenf}}\label{apevenf}
We prove a more general result: \begin{lemma} Assume that $n$ is even. Let $g \in C_0^\infty[0,\infty)$ and $$k_n(y,s) = \int\limits_0\limits^\infty \llg^{2n-3} N(s\llg) \int\limits_0\limits^\infty g(r) J(r \llg) dr d\llg.$$ Then $$k_n(y,s) =  \frac{(-1)^{\frac{n-2}{2}}}{\pi}  \int\limits_0\limits^\infty \left(\frac{d}{dr} \frac{1}{r}\right)^{n-1} g(y,r) \ln |r^2-s^2| dr.$$ \end{lemma}
Obviously, applying this lemma for $g(s)=(\mR_S f)(y,s)$ we obtain \mref{evenf}.

\bproof Since $n$ is even, due to \mref{wat}, we have $$N(s)= (-1)^{\frac{n-2}{2}}\left(\frac{1}{s}\frac{d}{ds} \right)^{\frac{n-2}{2}}N_0(s), ~J(s)= (-1)^{\frac{n-2}{2}}\left(\frac{1}{s}\frac{d}{ds} \right)^{\frac{n-2}{2}}J_0(s).$$
Thus, \begin{eqnarray*} k_n(y,s) = \int\limits_0\limits^\infty \llg^{2n-3} \left[\left(\frac{1}{s}\frac{d }{ds}\right)^{\frac{n-2}{2}}N_0 \right](s\llg) \int\limits_0\limits^\infty g(y,r) \left[\left(\frac{1}{s}\frac{d}{ds} \right)^{\frac{n-2}{2}}J_0 \right](r\llg) dr d\llg.\end{eqnarray*}
Since $$\left[\left(\frac{1}{s}\frac{d }{ds}\right)^{\frac{n-2}{2}}N_0 \right](s\llg)= \llg^{-(n-2)} \left(\frac{1}{s}\frac{d }{ds}\right)^{\frac{n-2}{2}} N_0(s\llg),$$ and $$\left[\left(\frac{1}{s}\frac{d }{ds}\right)^{\frac{n-2}{2}}J_0 \right](s\llg)= \llg^{-(n-2)} \left(\frac{1}{s}\frac{d }{ds}\right)^{\frac{n-2}{2}} J_0(s\llg),$$ we have

\begin{eqnarray*} k_n(y,s) &=&  \int\limits_0\limits^\infty \llg \left(\frac{1}{s}\frac{d }{ds}\right)^{\frac{n-2}{2}} N_0(s\llg) \int\limits_0\limits^\infty g(y,r) \left(\frac{1}{r}\frac{d}{dr} \right)^{\frac{n-2}{2}} J_0(r\llg) dr d\llg \\ &=& \left(\frac{1}{s}\frac{d }{ds}\right)^{\frac{n-2}{2}} \int\limits_0\limits^\infty \llg N_0(s\llg) \int\limits_0\limits^\infty g(y,r) \left(\frac{1}{r}\frac{d}{dr} \right)^{\frac{n-2}{2}} J_0(r\llg) dr d\llg.\end{eqnarray*}
Integrating by parts for the inner integral, we obtain \begin{eqnarray}\label{kn} k_n(y,s)= (-1)^{\frac{n-2}{2}}\left(\frac{1}{s}\frac{d }{ds}\right)^{\frac{n-2}{2}} l_n(y,s),\end{eqnarray} where \begin{eqnarray*} l_n(y,s)=\int\limits_0\limits^\infty \llg N_0(s\llg) \int\limits_0\limits^\infty \left(\frac{d}{dr} \frac{1}{r}\right)^{\frac{n-2}{2}}g(y,r)J_0(r\llg) dr d\llg.\end{eqnarray*}
Following the argument in \cite{AKKun}, we obtain \begin{eqnarray*} l_n(y,s) &=& - \frac{2}{\pi}  \int\limits_0\limits^\infty \left(\frac{d}{dr} \frac{1}{r}\right)^{\frac{n-2}{2}}g(y,r) \frac{1}{r^2-s^2} dr.\end{eqnarray*} Here the integral is understood in the principal value sense.
We now prove that \begin{eqnarray*} \left(\frac{1}{s}\frac{d}{ds}\right)^{\ag} l_n(y,s) = - \frac{2}{\pi}  \int\limits_0\limits^\infty \left(\frac{d}{dr} \frac{1}{r}\right)^{\frac{n-2}{2}+\ag}g(y,r) \frac{1}{r^2-s^2} dr,\end{eqnarray*} for any nonnegative integer $\ag$.
Indeed, we can rewrite \begin{eqnarray*} l_n(y,s) &=& - \frac{2}{\pi} \int\limits_0\limits^\infty \frac{1}{r} \left(\frac{d}{dr} \frac{1}{r}\right)^{\frac{n-2}{2}}g(y,r) \frac{r}{r^2-s^2} dr.\end{eqnarray*}
Integrating by parts (which can be justified for the principal value integral in question), we get \begin{eqnarray*} l_n(y,s) &=&  \frac{1}{\pi} \int\limits_0\limits^\infty \frac{d}{dr}\frac{1}{r} \left(\frac{d}{dr} \frac{1}{r}\right)^{\frac{n-2}{2}}g(y,r) \ln|r^2-s^2| dr.\end{eqnarray*}
Hence, \begin{eqnarray*} \left(\frac{1}{s} \frac{d}{ds}\right)l_n(y,s) &=&  \frac{2}{\pi} \int\limits_0\limits^\infty \frac{d}{dr}\frac{1}{r} \left(\frac{d}{dr} \frac{1}{r}\right)^{\frac{n-2}{2}}g(y,r) \frac{1}{s^2-r^2} dr \\ &=& -\frac{2}{\pi} \int\limits_0\limits^\infty \left(\frac{d}{dr} \frac{1}{r}\right)^{\frac{n-2}{2}+1}g(y,r) \frac{1}{r^2-s^2} dr.\end{eqnarray*}
By induction on $\ag$, we conclude that \begin{eqnarray*} \left(\frac{1}{s}\frac{d}{ds}\right)^{\ag} l_n(y,s) &=& - \frac{2}{\pi}  \int\limits_0\limits^\infty \left(\frac{d}{dr} \frac{1}{r}\right)^{\frac{n-2}{2}+\ag}g(y,r) \frac{1}{r^2-s^2} dr,\end{eqnarray*} for any nonnegative integer $\ag$. Therefore, due to \mref{kn}, \begin{eqnarray*}k_n(y,s) = \frac{2(-1)^{\frac{n}{2}}}{\pi}  \int\limits_0\limits^\infty \left(\frac{d}{dr} \frac{1}{r}\right)^{n-2}g(y,r) \frac{1}{r^2-s^2} dr \end{eqnarray*}
Integrating by parts again, we obtain \begin{eqnarray*}k_n(y,s) = \frac{(-1)^{\frac{n-2}{2}}}{\pi}  \int\limits_0\limits^\infty \left(\frac{d}{dr} \frac{1}{r}\right)^{n-1}g(y,r) \ln |r^2-s^2| dr.\end{eqnarray*} \eproof

\subsection{Proof of Lemma \ref{toprk}}
In view of (\ref{connect}), we have \begin{eqnarray}\label{doub1}\int\limits_0\limits^\infty \llg^{2n-3} N(s\llg) \int\limits_0\limits^\infty h(r) J(r \llg) dr d\llg &=& c \int\limits_0\limits^\infty \llg \mR(s,\llg) \int\limits_0\limits^\infty h(r) \mI(r, \llg) dr d\llg \\ \label{doub2} \int\limits_0\limits^\infty \llg^{2n-3} J(s\llg) \int\limits_0\limits^\infty h(r) N(r \llg) dr d\llg &=& c\int\limits_0\limits^\infty \llg \mI(s,\llg) \int\limits_0\limits^\infty h(r) \mR(r, \llg) dr d\llg.\end{eqnarray} Here the constant $c$ is the same in these two formulas.
Recall that, \begin{eqnarray*} G(s,\llg) = \mW(e_\llg)(s),\end{eqnarray*} where $e_{\llg}(s) =e^{i \llg s}$. Let $\cos_\llg(s)=\cos(\llg s)$ and $\sin_\llg(s)=\sin(\llg s)$, we then have
\begin{eqnarray*}\int\limits_0\limits^\infty \llg \mR(s,\llg) \int\limits_0\limits^\infty h(r) \mI(r, \llg) dr d\llg = \int\limits_0\limits^\infty \llg \mW(\cos_\llg)(s)\int\limits_0\limits^\infty h(r) \mW(\sin_\llg)(r) dr d \llg.\end{eqnarray*}
Similar to the argument in Lemma \ref{finchlem}, we can take $\mW$ out of the integral sign to obtain \begin{eqnarray}\label{bh}\int\limits_0\limits^\infty \llg \mR(s,\llg) \int\limits_0\limits^\infty h(r) \mI(r,\llg) dr d\llg = \mW(H)(s),\end{eqnarray} where \begin{eqnarray*} H(s) &=& \int\limits_0\limits^\infty \llg \cos (\llg s) \int\limits_0\limits^\infty h(r)\mW(\sin_\llg) dr d\llg \\ &=&  \frac{d}{ds}\int\limits_0\limits^\infty \sin (\llg s) \int\limits_0\limits^\infty h(r)\mW(\sin_\llg) dr d\llg.\end{eqnarray*}
Let $\mW^*$ be the $L^2-$adjoint of $\mW$, which is \begin{eqnarray}\label{mws} \mW^*(v)(s) =\left\{\begin{array}{l} (-1)^{\frac{n-2}{2}} c_n \int\limits_{0}^s \frac{\left(\frac{1}{t} \frac{d}{dt} \right)^{\frac{n-2}{2}} v(t)}{\sqrt{s^2-t^2}}  dt, \mbox{ if $n$ is even}, \\ (-1)^{\frac{n-3}{2}}c_n \left(\frac{1}{s} \frac{d}{ds} \right)^{\frac{n-3}{2}} \left(\frac{v(s)}{s}\right), \mbox{ if $n$ is odd}.
\end{array} \right.\end{eqnarray}
Given the condition in \reflemm{toprk}, integrating by parts (and changing the order of integration if $n$ is even), we obtain $$\int\limits_0\limits^\infty h(r)\mW(\sin_\llg) dr d\llg = \int\limits_0\limits^\infty \mW^*(h)(r)\sin_\llg(r) dr d\llg.$$
Therefore, \begin{eqnarray*}H(s) = \frac{d}{ds}\int\limits_0\limits^\infty \sin (\llg s) \int\limits_0\limits^\infty \mW^*(h)(r)\sin(\llg r) dr d\llg.\end{eqnarray*}
Since the inversion of the Fourier-sine transform is itself, we arrive at \begin{eqnarray*} H(s) = \frac{d}{ds}\mW^*(h)(s).\end{eqnarray*}
Thus, due to (\ref{bh}), \begin{eqnarray}\label{left}\int\limits_0\limits^\infty \llg \mR(s,\llg) \int\limits_0\limits^\infty h(r) J(r, \llg) dr d\llg = \left(\mW \frac{d}{ds} \mW^*\right)(h)(s).\end{eqnarray}

Similarly, since \begin{eqnarray*}\int\limits_0\limits^\infty \llg \mI(s,\llg) \int\limits_0\limits^\infty h(r) \mR(r, \llg) dr d\llg = \int\limits_0\limits^\infty \llg \mW(\sin_\llg)(s)\int\limits_0\limits^\infty h(r) \mW(\cos_\llg)(r) dr d \llg,\end{eqnarray*} we obtain
\begin{eqnarray*}&&\int\limits_0\limits^\infty \llg \mI(s,\llg) \int\limits_0\limits^\infty h(r) \mR(r, \llg) dr d\llg = \mW(H_1)(s),\end{eqnarray*} where \begin{eqnarray*}H_1(s) &=&\int\limits_0\limits^\infty \llg \sin (\llg s) \int\limits_0\limits^\infty \mW^*(h)(r)\cos(\llg r) dr d\llg \\ &=& - \frac{d}{ds}\int\limits_0\limits^\infty \cos (\llg s) \int\limits_0\limits^\infty \mW^*(h)(r)\cos(\llg r) dr d\llg \\ &=& -\frac{d}{ds}\mW^*(h)(s),\end{eqnarray*}
Here we have used the fact that the inversion of the Fourier-cosine transform is itself. Thus, \begin{eqnarray}\label{right} \int\limits_0\limits^\infty \llg \mR(s,\llg) \int\limits_0\limits^\infty h(r) J(r, \llg) dr d\llg = -\left(\mW \frac{d}{ds} \mW^*\right)(h)(s).\end{eqnarray}

From (\ref{doub1}), (\ref{doub2}), (\ref{left}), and (\ref{right}), we have \begin{eqnarray*}\int\limits_0\limits^\infty \llg^{2n-3} N(s\llg) \int\limits_0\limits^\infty h(r) J(r \llg) dr d\llg =- \int\limits_0\limits^\infty \llg^{2n-3} J(s\llg) \int\limits_0\limits^\infty h(r) N(r \llg) dr d\llg.\end{eqnarray*}

\reflemm{toprk} is proved.

\section{Remarks} \label{remm}
\begin{itemize}

\item The operator $\mW$ is an intertwining operator between the second derivative and Bessel operator (e.g., \cite{Tri,Lev,AKQ}), i.e. $$\mW\left[\left( \frac{d}{ds} \right)^2 v \right]= \left[\left( \frac{d}{ds}\right)^2 + \frac{n-1}{s} \frac{d}{ds}\right] \mW(v).$$  Thus, $\mW$ transforms the wave equation into the Darboux equation, which is known to describe spherical means (see \cite{CH2,Hel,As,Jo}).

\item The symmetry of $K(x,z,\llg)$ is similar to that of $I(x,z,\llg)$ in \cite{Kun07}. It implies Theorem \ref{main} as shown in Section \ref{direct}. It is not clear that the converse is true. However, a close look at the argument in Section \ref{direct} shows that Theorem \ref{main} implies the following weaker symmetry: $$\int\limits_{\RR}K(x,z,\llg) d\llg= \int\limits_{\RR}K(z,x,\llg) d\llg, \mbox{ for all } x,z \in B.$$

\item In \cite{FR07}, the authors derived an inversion formula using the known Neumann rather than Dirichlet observation data of the solution for wave solution for $n=3$ (see \cite[Theorem 12]{FR07}). Applying the results in this paper, one can derive the inversion formulas from Neumann data for any $n$. Indeed, looking at (\ref{mainap}) and (\ref{gradin2}), we arrive at: \begin{eqnarray*} f(x) = 2 \int\limits_{S}  \int\limits_{\RR} \overline{G}(|y-x|,\llg)\frac{ \pdh \hat{u}_0}{\pdh \nu_y}(y,\llg) d\llg d\sg(y). \end{eqnarray*} Due to Lemma \ref{finchlem}, we obtain \begin{eqnarray*} f(x) = 2 \int\limits_{S}  \mW(\pdh_{\nu_y}u)(y,|x-y|) d\sg(y), \end{eqnarray*} where $\pdh_{\nu_y}$ stands for the (outward pointing) normal derivative.
\end{itemize}
\section*{Acknowledgments} This work is partially supported by the NSF DMS grants 0604778 and 0715090, and the grant KUS-CI-016-04 from King Abdullah University of Science and Technology (KAUST). The author thanks the NSF and KAUST for the support.

The author is also grateful to Professor P. Kuchment for his suggestions and thanks Professor L. Kunyansky for useful discussions.

%\bibliographystyle{acm}
%\bibliography{research}

\medskip
 {\it E-mail address: }lnguyen@math.tamu.edu

\end{document}